\documentclass[german,a4paper,12pt]{article}
\usepackage{amsmath,amsthm,amsfonts,latexsym,amscd}

\usepackage{amsthm}

\usepackage{amssymb}
\swapnumbers
\theoremstyle{plain}

\theoremstyle{definition}

\theoremstyle{remark}

\begin{document}

{\bf Line bundles, connections, Deligne-Beilinson and absolute Hodge cohomology}\\

{\bf Helmut A. Hamm (M\"unster)}\\

\footnotesize
{\bf Abstract:} It is known that the Picard group of a complex manifold can be expressed as a Deligne cohomology group. One may wonder if the same holds for the Picard group of a smooth algebraic variety and Deligne-Beilinson cohomology but this is not true, as already remarked by M. Saito. We explain how one has to modify the latter, show that the Picard group can be expressed by absolute Hodge cohomology, too, and introduce an intermediate object between Picard group and usual Deligne-Beilinson cohomology group.\\
Similarly as in the case of Deligne cohomology one can relate line bundles with a regular connection to (modified) Deligne-Beilinson cohomology. In order to take irregular connections into account one has to change the definition of Deligne-Beilinson cohomology even more.\\

\normalsize

In this paper we are mainly interested in smooth complex algebraic varieties but start with complex manifolds.\\

It is well-known that Deligne cohomology of complex manifolds is closely related to line bundles and connections. \\
Let $M$ be a complex manifold. The Deligne cohomology group $H^k_{\cal D}(M,\mathbb{Z}(p))$ is defined to be the $k$-th hypercohomology group of the complex $\mathbb{Z}(p)\to \Omega^0_M\to\ldots\to \Omega^{p-1}_M\to 0\to\ldots$, where $\mathbb{Z}(p):=(2\pi i)^p\mathbb{Z}$, see [B1], [EV].\\
On the other hand, let $Pic\,M$ be the Picard group of $M$, i.e. the group of isomorphism classes of line bundles on $M$, $Pic_c\,M$ and $Pic_{ci}\,M$
the group of isomorphism classes of line bundles on $M$ with a connection resp. an integrable connection, see [HL2] for this notation. Then:\\

{\bf Theorem 0.1} (see [B1], [EV], [Ga]):\\ 
a) $H^2_{\cal D}(M,\mathbb{Z}(1))\simeq Pic\,M$,\\
b) $H^2_{\cal D}(M,\mathbb{Z}(2))\simeq Pic_c\,M$,\\
c) $H^2_{\cal D}(M,\mathbb{Z}(p))\simeq Pic_{ci}\,M, p\ge 3$.\\

We will study the question whether we can pass to the algebraic case, using Deligne-Beilinson cohomology $H^2_{DB}(X,\mathbb{Z}(p))$ instead of Deligne cohomology. 
It is defined as follows, see [B1], [EV]:\\

Let $X$ be a smooth complex algebraic variety; it is compactifiable according to Nagata [N]. Using resolution of singularities one finds a good compactification $\overline{X}$ of $X$, which means that $\overline{X}$ is smooth, too, and $D:=\overline{X}\setminus X$ is a divisor with normal crossings and smooth irreducible components. Then $H^ k_{DB}(X,\mathbb{Z}(p)):=\mathbb{H}^k(X^{an},\mathbb{Z}(p)_{DB})$, where $\mathbb{Z}(p)_{DB}:= cone(Rj_*^{an}\mathbb{Z}(p)\oplus F^p\Omega^\cdot_{\overline{X}^{an}}(log\,D)\to Rj_*^{an}\Omega^\cdot_{X^{an}})[-1]$; see [EV] Def. 2.6. Here $F^p$ denotes the Hodge filtration. For the use of cones in homological algebra see e.g. [GM]. Note that $H^ k_{DB}(X,\mathbb{Z}(p))$ is independent of the choice of the compactification, see [EV] Lemma 2.8.\\

We prefer a different description which is closer to the definition of Deligne cohomology:\\

{\bf Lemma 0.2} (see [S] p. 292 below): \\
$H^k_{DB}(X,\mathbb{Z}(p))=\mathbb{H}^k(\overline{X}^{an},cone(Rj^{an}_*\mathbb{Z}(p)\to \Omega^{\le p-1}_{\overline{X}^{an}}(log\,D))[-1])$.\\

{\bf Proof:} According to [EV] 2.7, $\mathbb{Z}(p)_{DB}$ is quasiisomorphic to $cone(Rj_*^{an}\mathbb{Z}(p)\to cone(F^p\Omega^\cdot_{\overline{X}^{an}}(log\,D)\to Rj_*^{an}\Omega^\cdot_{X^{an}}))[-1]$.\\
Here we may replace $Rj_*^{an}\Omega^\cdot_{X^{an}}$ by the quasiisomorphic complex $\Omega^\cdot_{\overline{X}^{an}}(log\,D)$. Finally $cone(F^p\Omega^\cdot_{\overline{X}^{an}}(log\,D)\to \Omega^\cdot_{\overline{X}^{an}}(log\,D))$ may be replaced by the quasiisomorphic complex $\Omega^{\le p-1}_{\overline{X}^{an}}(log\,D)$.\\

One might be optimistic because of the following reason: If $M$ is a complex manifold we have $Pic\,M\simeq H^1(M,{\cal O}^*_M)$. Now $H^0(M,{\cal O}^*_M)\simeq H^1_{\cal D}(M,\mathbb{Z}(1))$. For an algebraic variety $X$ we have similarly\\

{\bf Lemma 0.3} (see [EV] Prop. 2.12 (iii)): $H^1_{DB}(X,\mathbb{Z}(1))\simeq H^0(X,{\cal O}^*_X)$.\\

But there are difficulties with ${\cal O}^*_X$: it is not a coherent algebraic sheaf, and in the algebraic context there is no exponential sequence.\\

In fact M. Saito [S] has pointed out that $Pic\,X$ and $H^2_{DB}(X,\mathbb{Z}(1))$ need not be isomorphic, see [S] Remark 3.5 (i). I became aware of this paper by a hint of H. Esnault. See also the examples in \S 4.\\

As explained in [S] Prop. 3.4 the natural map $Pic\,X\to H^2_{DB}(X,\mathbb{Z}(1))$ is at least injective, see Lemma 2.14 below. We will give an explicit description of the cokernel.\\

Let us start with the following considerations. If $X$ is compact the Deligne-Beilinson cohomology of $X$ coincides with the Deligne cohomology of $X^{an}$, because of Lemma 0.2, so for $p=1$ we have the Picard group of $X$, for $p\ge 2$ the Picard group of line bundles with an integrable connection (note that in this case connections on line bundles are automatically integrable because of Hodge theory, see [HL2] Lemma 3.3). \\

Now let us return to the general (i.e. not necessarily compact) case. Each line bundle on $X$ corresponds to a Weil divisor, so it can be extended to $\overline{X}$, which means that $Pic\,\overline{X}\to Pic\,X$ is surjective, cf. [Ha] II Prop. 6.5. (This implies that the Chern class of a line bundle on $X$ comes from an element of $H^2(\overline{X}^{an};\mathbb{Z})$).\\

On the other hand we have an exact sequence
$$H^2(\overline{X}^{an},X^{an};\mathbb{Z})\to H^2_{DB}(\overline{X},\mathbb{Z}(1))\to H^2_{DB}(X,\mathbb{Z}(1))\to H^3(\overline{X}^{an},X^{an};\mathbb{Z})$$
see proof of Theorem 2.4 below.\\
We are no longer sure to have surjectivity for the middle arrow, so it is safer to replace $H^2_{DB}(X,\mathbb{Z}(1))$ by the kernel $H^2_{db}(X,\mathbb{Z}(1))$ of $H^2_{DB}(X,\mathbb{Z}(1))\to H^3(\overline{X}^{an},X^{an};\mathbb{Z})$. This will be done in the next section, in order to prove independence of the compactification it will be important that we may replace $H^3(\overline{X}^{an},X^{an};\mathbb{Z})$ by $H^3(\overline{X}^{an},X^{an};\mathbb{Q})$ here (cf. Lemma 1.5).\\

One of the main results will be (Theorem 2.4): 
$$Pic\,X\simeq H^2_{db}(X,\mathbb{Z}(1))$$ 
Instead of Deligne-Beilinson cohomology we can also look at absolute Hodge cohomology which is defined as follows:\\
$$H^2_{AH}(X,\mathbb{Z}(1)):=Ext^2_{MHS(\mathbb{Z})}(\mathbb{Z}^H,R\Gamma(X^{an},\mathbb{Z})(1))$$
where the index $H$ denotes Hodge and $MHS(\mathbb{Z})$ mixed Hodge structure over $\mathbb{Z}$.\\
Similarly, 
$$H^2_{p-AH}(X,\mathbb{Z}(1)):=Ext^2_{MHS(\mathbb{Z})^p}(\mathbb{Z}^H,R\Gamma(X^{an},\mathbb{Z})(1))$$
where the index $p$ denotes "polarized".\\
See [B2].\\

Now another important result will be:\\

{\bf Theorem 0.4:} $Pic\,X\simeq H^2_{p-AH}(X,\mathbb{Z}(1))\simeq H^2_{AH}(X,\mathbb{Z}(1))$.\\

{\bf 1. Essential Deligne-Beilinson cohomology}\\

Let $X$ be defined as in the previous section. As we have explained it is good to look at a modified version of Deligne-Beilinson cohomology:\\

{\bf Definition 1.1:} $H^k_{db}(X,\mathbb{Z}(p)):=ker(H^k_{DB}(X,\mathbb{Z}(p))\to H^{k+1}(\overline{X}^{an},X^{an};\mathbb{Q}))$ (essential Deligne-Beilinson cohomology).\\

Note that the mapping above factorizes through $H^k(X^{an};\mathbb{Q})$, and recall that $im(H^k(\overline{X}^{an};\mathbb{Q})\to H^k(X^{an};\mathbb{Q}))=W_kH^k(X^{an};\mathbb{Q})$.  See e.g. [E] Cor. 3.7.12. Therefore 
$$H^k_{db}(X,\mathbb{Z}(p))=ker(H^k_{DB}(X,\mathbb{Z}(p))\to H^k(X^{an};\mathbb{Q})/W_kH^k(X^{an};\mathbb{Q}))$$ 
This shows that $H^k_{db}(X,\mathbb{Z}(p))$ is independent of the compactification.\\

One has a more natural alternative to define $H^k_{db}(X,\mathbb{Z}(p))$:\\
$H^k_{db}(X,\mathbb{Z}(p)):=ker(H^k_{DB}(X,\mathbb{Z}(p))\to H^k(\overline{X}^{an},X^{an};\mathbb{Z}))$.\\
For $k=2$ this gives the same as before because of Corollary 1.6 below but in general there is no reason for coincidence.\\
The alternative definition does not depend on the choice of the compactification $\overline{X}$, too:\\
Let $\overline{X}_1, \overline{X}_2$ be two ``good'' compactifications of $X$, then there is a ``good'' compactification $\overline{X}_3$ and morphisms $f_i:\overline{X}_3\to \overline{X}_i, i=1,2$, which extend $id_X$: consider the closure $Z$ of the graph of $id_X$ in $\overline{X}_1\times\overline{X}_2$ and let $\overline{X}_3$ be obtained from $Z$ by resolution of singularities.\\
Therefore it is sufficient to show the following:\\
Let $f:\overline{X}_1\to\overline{X}_2$ be a morphism between two ``good'' compactifications of $X$ which extends $id_X$. Then $im(H^k(\overline{X}_i^{an};\mathbb{Z})\to H^k(X^{an};\mathbb{Z}))$ does not depend on $i, i=1,2$.\\
Now we apply [Do] Prop. 10.9, p. 311, with $h=f^{an}$, $K:=\overline{X}_2^{an}$, $\tilde{K}:=\overline{X}_2^{an}\setminus U$, $U$ being a good open neighbourhood of $\overline{X}_2^{an}\setminus X^{an}$ in $\overline{X}_2^{an}$ so that $H^k(\overline{X}_2^{an}\setminus U)\simeq H^k(X^{an})$, $L:=\tilde{L}:=\emptyset$. Let $j_i^U:X^{an}\setminus U\to \overline{X}_i^{an}$ be the inclusion, $i=1,2$.  We work with integral cohomology. Then we have a transfer map $(f^{an})^!$ which makes the following diagram commutative:\\
$$\begin{array}{ccc}
H^k(\overline{X}_1^{an})&\stackrel{(f^{an})^!}{\to}&H^k(\overline{X}_2^{an})\\
\downarrow (j_1^U)^*&&\downarrow (j_2^U)^*\\
H^k(X^{an}\setminus U)&\stackrel{(f^{an})^!}{\to}&H^k(X^{an}\setminus U)
\end{array}$$
Now the lower horizontal map is the identity, because $f^{an}|X^{an}\setminus U:X^{an}\setminus U\to X^{an}\setminus U$ is the identity. Hence we get a commutative diagram
$$\begin{array}{ccc}
H^k(\overline{X}_1^{an})&\stackrel{(f^{an})^!}{\to}&H^k(\overline{X}_2^{an})\\
\downarrow (j_1^{an})^*&&\downarrow (j_2^{an})^*\\
H^k(X^{an})&\stackrel{id}{\to}&H^k(X^{an})
\end{array}$$
where $j_i:X\to \overline{X}_i$ is the inclusion.\\
This implies: $(j_1^{an})^*H^k(\overline{X}_1^{an})=(j_2^{an})^*((f^{an})^!H^k(\overline{X}_1^{an}))\subset (j_2^{an})^*H^k(\overline{X}_2^{an})$,\\
whereas $(j_2^{an})^*H^k(\overline{X}_2^{an})=(j_1^{an})^*((f^{an})^*H^k(\overline{X}_2^{an}))\subset (j_1^{an})^*H^k(\overline{X}_1^{an})$,\\
so $(j_1^{an})^*H^k(\overline{X}_1^{an})=(j_2^{an})^*H^k(\overline{X}_2^{an})$.\\

{\bf I.} First let us study the case $k=0$. \\

We have $H^0_{db}(X,\mathbb{Z}(p))=H^0_{DB}(X,\mathbb{Z}(p))$ because $H^1(\overline{X}^{an},X^{an};\mathbb{Q})=0$.\\

So we need only to compute $H^0_{DB}(X,\mathbb{Z}(p))$:\\

{\bf Lemma 1.2:} a) $H^0_{DB}(X,\mathbb{Z}(0))=H^0(X^{an};\mathbb{Z})$,\\
b) $H^0_{DB}(X,\mathbb{Z}(p))=0$ for $p>0$.\\

{\bf Proof:} a) obvious.\\
b) We have an exact sequemce
$$0\to H^0_{DB}(X,\mathbb{Z}(1))\to H^0_{DB}(X;\mathbb{Z}(0))\to H^0(\overline{X}^{an},{\cal O}_{\overline{X}^{an}})$$
where the second map coincides with the injective mapping $H^0(X^{an};\mathbb{Z})\simeq H^0(\overline{X}^{an},\mathbb{Z}_{\overline{X}^{an}})\to H^0(\overline{X}^{an},{\cal O}_{\overline{X}^{an}})$. So 
$H^0_{DB}(X,\mathbb{Z}(1))=0$.\\
Obviously, $H^0_{DB}(X,\mathbb{Z}(p))\simeq H^0_{DB}(X,\mathbb{Z}(1))$ for $p>1$.\\

{\bf II.} Consider now the case $k=1$: First we have:\\

{\bf Lemma 1.3:} a) $H^1_{DB}(X,\mathbb{Z}(0))=H^1(X^{an};\mathbb{Z})$.\\
b) (see Lemma 0.3) $H^1_{DB}(X,\mathbb{Z}(1))=H^0(X,{\cal O}^*_X)$.\\
c) $H^1_{DB}(X,\mathbb{Z}(p))=H^0(X^{an};\mathbb{C}^*)$ for $p>1$.\\

{\bf Proof:} a) obvious.\\
b) see loc.cit.\\
c) We have an exact sequence
$$0\to H^1_{DB}(X,\mathbb{Z}(2))\to H^1_{DB}(X,\mathbb{Z}(1))\to H^0(\overline{X},\Omega^1_{\overline{X}}(log\,D))$$
where the last map corresponds to $H^0(X,{\cal O}^*_X)\to H^0(\overline{X},\Omega^1_{\overline{X}}(log\,D))$: $g\mapsto \frac{dg}{g}$.\\
So $H^1_{DB}(X,\mathbb{Z}(2))\simeq H^0(X,\mathbb{C}^*_X)\simeq H^0(X^{an};\mathbb{C}^*)$.\\
For $p>2$ we have $H^1_{DB}(X,\mathbb{Z}(p))\simeq H^1_{DB}(X,\mathbb{Z}(2))$.\\

{\bf Proposition 1.4:} a) $H^1_{db}(X,\mathbb{Z}(0))=im(H^1(\overline{X}^{an};\mathbb{Z})\to H^1(X^{an};\mathbb{Z}))$.\\
b) $H^1_{db}(X,\mathbb{Z}(p))=H^0(X^{an};\mathbb{C}^*)$ for $p\ge 1$.\\

{\bf Proof:} a) Since $H^2(\overline{X}^{an},X^{an};\mathbb{Z})$ is torsion free, $H^1_{db}(X,\mathbb{Z}(0))$ is the kernel of $H^1_{DB}(X,\mathbb{Z}(0))\to H^2(\overline{X}^{an},X^{an};\mathbb{Z})$, i.e. of $H^1(X^{an};\mathbb{Z})\to H^2(\overline{X}^{an},X^{an};\mathbb{Z})$.\\
The rest is obvious.\\
b) $H^1_{db}(X,\mathbb{Z}(1))$ is the kernel of $H^1_{DB}(X,\mathbb{Z}(1))\to H^2(\overline{X}^{an},X^{an};\mathbb{C})$.\\
Now identify $H^1_{DB}(X,\mathbb{Z}(1))$ with $H^0(X,{\cal O}^*_X)$ and $H^2(\overline{X}^{an},X^{an};\mathbb{C})$ with \\
$Hom(H_2(\overline{X}^{an},X^{an};\mathbb{Z}),\mathbb{C})$. Then the mapping is given by $g\mapsto ([c]\mapsto \int_c \frac{dg}{g})$, where $c$ is a singular $1$-cycle in $X^{an}$ which is a boundary in $\overline{X}^{an}$.\\
In order to see this, note that we have a commutative diagram
$$\begin{array}{ccc}
H^0(X,{\cal O}^*_X)&&\\
\downarrow&\searrow&\\
H^0(X^{an},{\cal O}^*_{X^{an}})&\to&Hom(H_2(\overline{X}^{an},X^{an};\mathbb{Z}),\mathbb{C})
\end{array}$$
Now let $g\in H^0(X,{\cal O}^*_X)$. Then $g$ can be considered as a rational function on $\overline{X}$. If $g$ is in the kernel, its divisor whose support is contained in $D$ must be $0$, hence $g$ extends to an element of $H^0(\overline{X},{\cal O}^*_{\overline{X}})=H^0(\overline{X},\mathbb{C}^*_{\overline{X}})$. So $H^1_{db}(X,\mathbb{Z}(1))\simeq H^0(X^{an};\mathbb{C}^*)$.\\
The case $p\ge 2$ is clear: any element of $H^1_{DB}(X,\mathbb{Z}(p))\simeq H^0(X^{an};\mathbb{C}^*)$ is mapped onto $0\in H^1(X^{an};\mathbb{C})$, because $H^0(X^{an};\mathbb{C}^*)\to H^1(X^{an};\mathbb{Z})\to H^1(X^{an};\mathbb{C})$ is exact, hence onto $0\in H^2(\overline{X}^{an},X^{an};\mathbb{C})$.\\

{\bf III.} The interesting case is $k=2$. The fact that $H^2_{db}(X,\mathbb{Z}(p))$ is independent of the choice of the compactification will be confirmed by the relation to Picard groups. Now we have:\\

{\bf Lemma 1.5:} $H^3(\overline{X}^{an},X^{an};\mathbb{Z})$ is torsion free, i. e. $H^3(\overline{X}^{an},X^{an};\mathbb{Z})\to H^3(\overline{X}^{an},X^{an};\mathbb{C})$ is injective.\\

{\bf Proof:} Note that $H^3(\overline{X}^{an},X^{an};\mathbb{Z})\to H^3(\overline{X}^{an},X^{an};\mathbb{C})$ is injective if and only if $H^2(\overline{X}^{an},X^{an};\mathbb{C})\to H^2(\overline{X}^{an},X^{an};\mathbb{C}^*)$ is surjective, because of the exponential sequence. This is the case: look at the commutative diagram
$$\begin{array}{ccc}
H^2(\overline{X}^{an},X^{an};\mathbb{C})&\to& H^2(\overline{X}^{an},X^{an};\mathbb{C}^*)\\
\downarrow\simeq&&\downarrow\simeq\\
H_{2n-2}(D^{an};\mathbb{C})&\to &H_{2n-2}(D^{an};\mathbb{C}^*)\\
\downarrow\simeq&&\downarrow\simeq\\
\mathbb{C}^k&\to&(\mathbb{C}^*)^k
\end{array}$$
where $k$ is the number of irreducible components of $D$ and $n=\dim\,X$ (without loss of generality we may assume that $X$ is purely $n$-dimensional). Note that for singular cohomology, $H_{2n-2}(\tilde{D}^{an})\simeq H_{2n-2}(D^{an})$, where $\tilde{D}$ is the (non-singular) normalization of $D$ and the coefficients are arbitrary, because $H_{2n-2}(\tilde{D}^{an})\simeq H_{2n-2}(\tilde{D}^{an},\tilde{\Sigma}^{an})\simeq H_{2n-2}(D^{an},\Sigma^{an})\simeq H_{2n-2}(D^{an})$ where $\Sigma$ is the singular locus of $D$ and $\tilde{\Sigma}$ its inverse image in $\tilde{D}$.\\
The lowest horizontal arrow is surjective.\\
Alternative proof:\\
By [Sp] Cor. 5.5.4, p. 244, the torsion subgroups of $H^3(\overline{X}^{an},X^{an};\mathbb{Z})$ and $H_2(\overline{X}^{an},X^{an};\mathbb{Z})$ are isomorphic.\\
Now $H_2(\overline{X}^{an},X^{an};\mathbb{Z})\simeq H^{2n-2}(D^{an};\mathbb{Z})$, by [Sp] Theorem 6.2.17, p. 296, and Cor. 6.1.11, p. 291. Then $H^{2n-2}(D^{an};\mathbb{Z})\simeq H^{2n-2}(D^{an},\Sigma^{an};\mathbb{Z})$, because $H^{k}(\Sigma^{an};\mathbb{Z})=0$ for $k=2n-3,2n-2$. And $H^{2n-2}(D^{an},\Sigma^{an};\mathbb{Z})\simeq H_0(D^{an}\setminus\Sigma^{an};\mathbb{Z})$ is torsion free.\\

{\bf Corollary 1.6:} $H^2_{db}(X,\mathbb{Z}(p))=ker(H^2_{DB}(X,\mathbb{Z}(p))\to H^3(\overline{X}^{an},X^{an};\mathbb{Z}))$\\

In the next section we will show that $H^2_{db}(X,\mathbb{Z}(1))\simeq Pic\,X$.\\

{\bf 2. Line bundles on smooth algebraic varieties}\\

Let $X$ be a smooth complex algebraic variety (not necessarily irreducible). Choose a good compactification as before.\\

Before turning to the main result about the algebraic Picard group it is useful to have some preparation:\\

Let $j:X\to\overline{X}$ be the inclusion.\\

{\bf Lemma 2.1:} $Pic\,X\simeq H^1(\overline{X},j_*{\cal O}^*_X)$, and there is an exact sequence
$$0\to H^0(\overline{X},{\cal O}^ *_{\overline{X}})\to H^0(X,{\cal O}^ *_X)\to H^0(\overline{X},j_*{\cal O}_X/{\cal O}^ *_{\overline{X}})\to Pic\,\overline{X}\to Pic\,X\to 0$$
{\bf Proof:} Look at the exact cohomology sequence for
$$0\to{\cal O}^*_{\overline{X}}\to j_*{\cal O}^*_X\to j_*{\cal O}^*_X/{\cal O}^*_{\overline{X}}\to 0$$
First we have $H^0(X,{\cal O}^*_X)=H^0(\overline{X},j_*{\cal O}^*_X)$. \\
Now we have two alternatives:\\
a) We know already that the sequence
$$0\to H^0(\overline{X},{\cal O}^ *_{\overline{X}})\to H^0(X,{\cal O}^ *_X)\to H^0(\overline{X},j_*{\cal O}_X/{\cal O}^ *_{\overline{X}})\to Pic\,\overline{X}\to H^1(\overline{X},j_*{\cal O}^*_X)$$
is exact.\\
Let us first show that we may replace here $H^1(\overline{X},j_*{\cal O}^*_X)$ by $Pic\,X$, without destroying the exactness:\\
We have $Pic\,X=H^1(X,{\cal O}^*_X)=\mathbb{H}^1(\overline{X},Rj_*{\cal O}^*_X)$. There is a commutative diagram
$$\begin{array}{ccc} 
H^1(\overline{X},{\cal O}^*_{\overline{X}})&\to& H^1(\overline{X},j_*{\cal O}^*_X)\\
&\searrow&\downarrow\\
&&H^1(X,{\cal O}^*_X)
\end{array}$$
Now there is a spectral sequence $E^{pq}_2=H^p(X,R^qj_*{\cal O}^*_X)\Rightarrow H^{p+q}(X,{\cal O}^*_X)$, see [Go] II 4.5, p. 177, and $H^1(\overline{X},j_*{\cal O}^*_X)=E^{10}_2\to H^1(X,{\cal O}^*_X)$ is injective, see [Go] I Th. 4.5.1, p. 82. So we may indeed replace $H^1(\overline{X},j_*{\cal O}^*_X)$ by $Pic\,X$ in the sequence above.\\
Now it is well-known that $Pic\,\overline{X}\to Pic\,X$ is surjective, see [Ha] II Prop. 6.5. So we obtain the desired exact sequence; furthermore, $H^1(\overline{X},j_*{\cal O}_X^*)\to H^1(X,{\cal O}_X^*)\simeq Pic\,X$ must be bijective.

\vspace{2mm}
b) First, $R^1j_*{\cal O}^*_X=0$: Let $x\in\overline{X}, s_x\in (R^1j_*{\cal O}^*_X)_x$, represented by $s\in H^1(U,j_*{\cal O}^*_X)$. Let $s_1$ be the image of $s$ in $H^1(U\cap X,{\cal O}^*_X)\simeq Pic(U\cap X)$. Now $Pic\,U\to Pic\,U\cap X$ is surjective, hence $s_1$ comes from $\hat{s}_1\in Pic\,U$. Then $\hat{s}_1$ comes from a line bundle on $U$. On some neighbourhood of $x$ the latter is trivial, so after shrinking $U$ if necessary we may assume that $\hat{s}_1=0$, hence $s_1=0$. As in the proof of a), the mapping $H^1(U,j_*{\cal O}^*_X)\to H^1(U\cap X,{\cal O}^*_X)$ is injective, so $s=0$, hence $s_x=0$.\\
This implies that $Pic\,X=H^1(\overline{X},j_*{\cal O}^*_X)$.\\
Furthermore the sheaf $j_*{\cal O}^*_X/{\cal O}^*_{\overline{X}}$ is flabby, because we may assume without less of generality that $X$ is irreducible and because the sheaf under consideration is the sheaf of Cartier divisors with support in $D$; cf. [Go] II Example 3.1.1, p. 147. Therefore $H^1(\overline{X},j_*{\cal O}^*_X/{\cal O}^*_{\overline{X}})=0$.\\
We obtain now the desired exact sequence as part of the long exact cohomology sequence for 
$0\to{\cal O}^*_{\overline{X}}\to j_*{\cal O}^*_X\to j_*{\cal O}^*_X/{\cal O}^*_{\overline{X}}\to 0$.\\

{\bf Remark 2.2:} $H^k(X,{\cal O}^*_X)=0, k\ge 2$, and $R^kj_*{\cal O}^*_X=0, k\ge 1$.\\

In order to show this look at the long exact cohomology sequence for\\
$0\to {\cal O}^*_X\to {\cal M}^*_X\to {\cal M}^*_X/{\cal O}^*_X\to 0$.\\
The sheaves ${\cal M}^*_X$ and ${\cal M}^*_X/{\cal O}^*_X$ are flabby: the connected components of $X$ are irreducible because $X$ is smooth, so we may suppose without loss of generality that $X$ is irreducible. Then ${\cal M}_X$ is constant, see [Ha] II proof of Prop. 6.15, p. 145, hence ${\cal M}^*_X$, too. Furthermore note that ${\cal M}^*_X/{\cal O}^*_X$ is the sheaf of Cartier divisors. The rest follows from [Go] II Example 3.1.1, p. 147.\\
We know by part b) in the proof of Lemma 2.1 that $R^1j_*{\cal O}^*_X=0$. Furthermore, $R^kj_*{\cal M}^*_X=R^kj_*({\cal M}^*_X/{\cal O}^*_X)=0, k\ge 1$, because we deal with flabby sheaves, hence $R^kj_*{\cal O}^*_X=0, k\ge 2$.\\

{\bf Lemma 2.3:} The complex $\tilde{\mathbb{Z}}(1):=cone(F^1\Omega^{\cdot}_{\overline{X}^{an}}(log\,D)\to j_*^{an}cone(\mathbb{Z}(1)\to \Omega^\cdot_{X^{an}}))[-1]$ is a resolution of $j_*^m{\cal O}^*_{X^{an}}$.\\

{\bf Proof:} We proceed similarly as in the proof of [EV] Prop. 2.12:\\
First, $\tilde{\mathbb{Z}}(1)$ is quasiisomorphic to\\
$\Omega^1_{\overline{X}^{an}}(log\,D)\oplus j_*^{an}{\cal O}^*_{X^{an}}\stackrel{\Delta_0}{\to}  \Omega^2_{\overline{X}^{an}}(log\,D)\oplus j_*^{an}\Omega^1_{X^{an}}\stackrel{\Delta_1}{\to}
 \Omega^3_{\overline{X}^{an}}(log\,D)\oplus j_*^{an}\Omega^2_{X^{an}}\stackrel{\Delta_2}{\to}\ldots$\\
where $\Delta_0(\omega,f):=(d\omega,\omega-\frac{df}{f}), \Delta_i(\omega,\phi):=(d\omega,\omega-d\phi), i>0$, see [EV] loc.cit.\\
So $(\omega,f)\in ker\Delta_0$ iff $\omega=\frac{df}{f}$. Suppose that $\omega$ is a closed logarithmic $1$-form: $\omega=\sum_{j=1}^s c_j\frac{dz_j}{z_j}+\phi$, $\phi$ holomorphic. Then $d\phi=0$, so $\phi=dg$. Then $f$ must be of the form $cz_1^{c_1}\cdots z_s^{c_s}e^g$ with some $c\in\mathbb{C}^*$. Since $f$ is uni-valued we have that $c_j\in\mathbb{Z}$, i.e. $f$ is meromorphic along $D$. This means that $ker\Delta_0\simeq j_*^m{\cal O}^*_{X^{an}}$.\\
The cone of $\Omega^\cdot_{\overline{X}^{an}}(log\,D)\to j_*^{an}\Omega^\cdot_{X^{an}}$ is acyclic. By comparison we obtain that $ker\Delta_i=im\Delta_{i-1}, i\ge 2$. This holds for $i=1$, too: Suppose that $(\omega,\phi)\in (ker\,\Delta_1)_x$. Then there is a $(\theta,g)\in \Omega^1_{\overline{X}^{an},x}(log\,D)\oplus (j_*^{an}{\cal O}_{X^{an}})_x$ such that $(\omega,\phi)=(d\theta,\theta-dg)$, again because the cone of $\Omega^\cdot_{\overline{X}^{an}}(log\,D)\to j_*^{an}\Omega^\cdot_{X^{an}}$ is acyclic. Put $f:=e^g$. Then $\Delta_0(\theta,f)=(\omega,\phi)$.\\

{\bf Theorem 2.4:} $H^2_{db}(X,\mathbb{Z}(1))\simeq Pic\,X$.\\

{\bf Proof:} Note that $\Omega^0_{\overline{X}}(log\,D)\simeq {\cal O}_{\overline{X}}$. The morphism
$cone(\mathbb{Z}(1)_{\overline{X}^{an}}\to {\cal O}_{\overline{X}^{an}})[-1]\to cone(Rj^{an}_*\mathbb{Z}(1)_{X^{an}}\to {\cal O}_{\overline{X}^{an}})[-1])$ induces a homomorphism 
$H^2_{DB}(\overline{X},\mathbb{Z}(1))\to H^2_{DB}(X,\mathbb{Z}(1))$
which fits into an exact sequence
$$H^2(\overline{X}^{an},X^{an};\mathbb{Z})\to H^2_{DB}(\overline{X},\mathbb{Z}(1))\to H^2_{DB}(X,\mathbb{Z}(1))\to H^3(\overline{X}^{an},X^{an};\mathbb{Z})$$
Note that $H^2_{db}(\overline{X},\mathbb{Z}(p))=H^2_{DB}(\overline{X},\mathbb{Z}(p))$.\\
Obviously we obtain an exact sequence
$$H^2(\overline{X}^{an},X^{an};\mathbb{Z})\to H^2_{db}(\overline{X},\mathbb{Z}(1))\to H^2_{db}(X,\mathbb{Z}(1))\to 0$$
Now $Pic(\overline{X})\to Pic\,X$ is surjective, and we have an exact sequence
$$H^0(\overline{X},j_*{\cal O}^*_X/{\cal O}^*_{\overline{X}})\to Pic(\overline{X})\to Pic(X)\to 0$$
see Lemma 2.1. Altogether, $Pic(X)\simeq H^2_{db}(X,\mathbb{Z}(1))$ - as soon as the two exact sequences fit together to a commutative diagram
$$\begin{array}{ccccccc}
H^0(\overline{X},j_*{\cal O}^*_X/{\cal O}^*_{\overline{X}})&\to&Pic(\overline{X})&\to&Pic(X)&\to&0\\
\downarrow\simeq&&\downarrow\simeq&&\downarrow\\
H^2(\overline{X}^{an},X^{an};\mathbb{Z})&\to&H^2_{db}(\overline{X},\mathbb{Z}(1))&\to&H^2_{db}(X,\mathbb{Z}(1))&\to&0
\end{array}$$

Now let us show that this is the case:\\
We have a complex analytic analogue of the short exact sequence in Lemma 2.1:
$$0\to{\cal O}^*_{\overline{X}^{an} }\to j_*^m{\cal O}^*_{X^{an}}\to j^m_*{\cal O}^*_{X^{an}}/{\cal O}^*_{\overline{X}^{an}}\to 0$$
where $j_*^m\Omega^p_{X^{an}}$ is the sheaf of differential $p$-forms which are holomorphic an $X^{an}$ and meromorphic on $\overline{X}^{an}$.
This leads to a commutative diagram:
\footnotesize
$$\begin{array}{ccccc}
H^0(\overline{X},j_*{\cal O}_X^*/{\cal O}^ *_{\overline{X}})&\to&Pic\,\overline{X}&\to&Pic\,X\\
\downarrow\simeq&&\downarrow\simeq&&\downarrow\\
H^0(\overline{X}^{an},j_*^m{\cal O}^*_{X^{an}}/{\cal O}^*_{\overline{X}^{an}})&\to&Pic\,\overline{X}^{an}&\to&H^1(\overline{X}^{an},j_*^m{\cal O}^*_{X^{an}})
\end{array}$$
\normalsize
The second vertical arrow is bijective by GAGA. The first vertical is bijective, too, because in both cases we are dealing with Cartier divisors on $\overline{X}$ with support in $D$.\\
Now ${\cal O}^*_{\overline{X}^{an}}$ is quasiisomorphic to $cone(\mathbb{Z}(1)_{\overline{X}^{an}}\to{\cal O}_{\overline{X}^{an}})[-1])$, hence to $cone(F^1\Omega^{\cdot}_{\overline{X}^{an}}(log\,D)\to cone(\mathbb{Z}(1)_{\overline{X}^{an}}\to j_*^{an}\Omega^\cdot_{X^{an}}))[-1]$,\\
and $j_*^ m{\cal O}^*_{X^{an}}$ is quasiisomorphic to $\tilde{\mathbb{Z}}(1)$:= $cone(F^1\Omega^{\cdot}_{\overline{X}^{an}}(log\,D)\to j_*^{an}cone(\mathbb{Z}(1)\to \Omega^\cdot_{X^{an}}))[-1]$, see Lemma 2.3.\\
Therefore $j_*^m{\cal O}^*_{X^{an}}/{\cal O}^*_{\overline{X}^{an}}$ is quasiisomorphic to 
$cone(cone(\mathbb{Z}(1)\to j_*^{an}\Omega^\cdot_{X^{an}})[-1]\to j_*^{an}cone(\mathbb{Z}(1)\to \Omega^\cdot_{X^{an}})[-1])$,\\
hence to $cone(cone(\mathbb{Z}(1)\to j_*^{an}{\cal O}_{X^{an}})[-1]\to j_*^{an}{\cal O}^*_{X^{an}})$, i.e. to $R^1j_*^{an}\mathbb{Z}(1)$.\hfill{(*)}\\
Note here that $\mathbb{Z}(1)\simeq j_*^{an}\mathbb{Z}(1)$, and
$$0\to j_*^{an}\mathbb{Z}(1)\to j_*^{an}{\cal O}_{X^{an}}\to j_*^{an}{\cal O}^*_{X^{an}}\to R^1j_*^{an}\mathbb{Z}(1)\to 0$$
is exact.\\
Now $H^0(\overline{X}^{an}, R^1j_*^{an}\mathbb{Z}(1))\simeq H^0(\overline{X}^{an},R^2\Gamma_{D^{an}}\mathbb{Z}(1))\simeq H^2(\overline{X}^{an},X^{an};\mathbb{Z}(1))$, because $R^k\Gamma_{D^{an}}\mathbb{Z}(1))=0, k<2$.\\
Therefore $H^0(\overline{X}^{an},j_*^m{\cal O}^*_{X^{an}}/{\cal O}^*_{\overline{X}^{an}})\simeq H^2(\overline{X}^{an},X^{an};\mathbb{Z}(1))$.\\
Now map $\tilde{\mathbb{Z}}(1)$ to $\mathbb{Z}(1)_{DB}\simeq cone(F^1\Omega^{\cdot}_{\overline{X}^{an}}(log\,D)\to cone(Rj_*^{an}(\mathbb{Z}(1)\to \Omega^\cdot_{X^{an}}))[-1]$.\\
Then the quotient $j_*^m{\cal O}^*_{X^{an}}/{\cal O}^*_{\overline{X}^{an}}$ is mapped to $cone(\mathbb{Z}(1)\to Rj_*^{an}\mathbb{Z}(1))[-1]$.\\
The group $H^0$ of this complex is isomorphic to $H^1(\overline{X}^{an},cone(\mathbb{Z}(1)\to Rj_*^{an}\mathbb{Z}(1)))=\mathbb{H}^2(\overline{X}^{an},R\Gamma_{D^{an}}\mathbb{Z}(1))=H^2(\overline{X}^{an},X^{an};\mathbb{Z}(1))$.\\
Altogether we obtain a commutative diagram
\footnotesize
$$\begin{array}{ccccccc}
H^0(\overline{X},j_*{\cal O}_X/{\cal O}^ *_{\overline{X}})&\to&Pic\,\overline{X}&\to&Pic\,X&\to&0\\
\downarrow\simeq&&\downarrow\simeq&&\downarrow&&\downarrow\\
H^2(\overline{X}^{an},X^{an};\mathbb{Z}(1))&\to&H^2_D(\overline{X},\mathbb{Z}(1))&\to&H^2_{DB}(X,\mathbb{Z}(1))&\to&H^3(\overline{X}^{an},X^{an};\mathbb{Z})
\end{array}$$
\normalsize
where it is easy to verify that the square on the right is commutative, too. \\
So replace $H^2_{DB}(X,\mathbb{Z}(1))$ by $H^2_{db}(X,\mathbb{Z}(1)):=ker(H^2_{DB}(X,\mathbb{Z}(1))\to H^3(\overline{X}^{an},X^{an};\mathbb{Z})$, then we obtain our statement.\\

In order to derive Theorem 0.4 let us quote a result by M. Saito:\\

{\bf Proposition 2.5:} (see [S] Prop. 3.4, p. 294) The natural mappings
$$Pic\,X\to H^2_{p-AH}(X,\mathbb{Z}(1))\to H^2_{AH}(X,\mathbb{Z}(1))\to H^2_{DB}(X,\mathbb{Z}(1))$$
are all injective, the middle map is bijective, and the first map has a finite cokernel.\\

{\bf Proof of Theorem 0.4:} We use Proposition 2.5.\\
In particular, there is an injective mapping $H^2_{AH}(X,\mathbb{Z}(1))\to H^2_{DB}(X,\mathbb{Z}(1))$.\\
Altogether, $H^2_{AH}(X,\mathbb{Z}(1))/Pic\,X\subset H^2_{DB}(X,\mathbb{Z}(1))/Pic\,X\subset H^3(\overline{X}^{an},X^{an};\mathbb{Z})$, hence the left hand group is finite and torsion free, by Lemma 1.5, hence $0$, so we have $Pic\,X\simeq H^2_{AH}(X,\mathbb{Z}(1))$.\\

{\bf Remark:} In [S] Remarks 3.5 (i), p. 297, we find comments about the injections:

\vskip.1in
a) It is said that the cokernel of $Pic\,X\to H^2_{AH}(X,\mathbb{Z}(1))$ is non-trivial as soon as the group $H^2(X^{an},\mathbb{Z})/W_2H^2(X^{an},\mathbb{Z})$ is not torsion free. But the latter cannot occur: Note that $W_2H^2(X^{an};\mathbb{Z})$ is in [S] loc. cit. the image of $H^2(\overline{X}^{an};\mathbb{Z})\to H^2(X^{an};\mathbb{Z})$, as we will see. But $H^2(X^{an},\mathbb{Z})/W_2H^2(X^{an},\mathbb{Z})$ is then isomorphic to a subgroup of $H^3(\overline{X}^{an},X^{an};\mathbb{Z})$ which is torsion free, by Lemma 1.5.\\
That $W_2H^2(X^{an};\mathbb{Z})$ is the image of $H^2(\overline{X}^{an};\mathbb{Z})\to H^2(X^{an};\mathbb{Z})$ can be seen as follows: From the definition of the weight filtration it is derived that $Gr_2^WH^2(X^{an};\mathbb{Z})=coker(\oplus H^0(D_i^{an};\mathbb{Z})\to H^2(\overline{X}^{an};\mathbb{Z}))$, and \\
$Gr_2^WH^i(X^{an};\mathbb{Z})=0, i<2$; see [S] (3.4.2), p. 295. \\
Therefore
$W_2H^2(X^{an};\mathbb{Z})=Gr_2^WH^2(X^{an};\mathbb{Z})\\
=coker(H^2(\overline{X}^{an},X^{an};\mathbb{Z})\to H^2(\overline{X}^{an};\mathbb{Z}))\\
=H^2(\overline{X}^{an};\mathbb{Z})/ker(H^2(\overline{X}^{an};\mathbb{Z})\to H^2(X^{an};\mathbb{Z}))\\
=im(H^2(\overline{X}^{an};\mathbb{Z})\to H^2(X^{an};\mathbb{Z}))$.

\vskip.1in
b) It is said that the cokernel of the mapping $H^2_{AH}(X,\mathbb{Z}(1))\to H^2_{DB}(X,\mathbb{Z}(1))$ is not torsion (i.e. not a finite group) under some Hodge theoretic condition, in particular if $X$ is the complement of an elliptic curve in $\mathbb{P}_2$: See the comment after Lemma 2.14 below; our example b) in section 4 corresponds to Saito's example.\\

Let us turn to related results:\\

Note that $H^k(X,{\cal O}_X)\simeq H^k(\overline{X},j_*{\cal O}_X)$, because $X\mapsto \overline{X}$ is affine, and $H^k(\overline{X},j_*{\cal O}_X)\simeq H^k(\overline{X}^{an},j_*^m{\cal O}_{X^{an}})$ by GAGA. Similarly, $H^k(X^{an},{\cal O}_{X^{an}})\simeq H^k(\overline{X}^{an},j_*^{an}{\cal O}_{X^{an}})$ because $j^{an}$ is Stein.\\
The situation is completely different if we pass from ${\cal O}$ to ${\cal O}^*$. Anyhow we have $R^kj_*{\cal O}^*_X=0$, $k>0$, see Remark 2.2, and $H^k(X,{\cal O}^*_X)\simeq H^k(\overline{X},j_*{\cal O}^*_X)=0$ for all $k\ge 2$. Furthermore, for $k=0$ we still have $H^0(\overline{X},j_*{\cal O}^*_X)\simeq H^0(\overline{X}^{an},j_*^m{\cal O}^*_{X^{an}})$. But the situation for $k=1$ is much more complicated:\\
If $X$ is compact we still have $H^k(X,{\cal O}^*_X)\simeq H^k(X^{an},{\cal O}^*_{X^{an}})$ for $k=1$ by GAGA. But if $X$ is non-compact this no longer holds, see Example a),b).\\ 
The same holds if $X$ is compact but $k\ge 2$, see Example c) in section 4.\\

{\bf Proposition 2.6:} a) Let $D=D_1\cup\ldots\cup D_r$ be the decomposition of $D$ into irreducible components. Then we have an exact sequence
$$0\to Pic\,X\to H^1(\overline{X}^{an},j_*^m{\cal O}^*_{X^{an}})\to \oplus_{j=1}^rH^1(D_i^{an};\mathbb{Z})\to H^2(\overline{X}^{an},{\cal O}^*_{\overline{X}^{an}})$$
b) There is an exact sequence
$$0\to Pic\,X\to H^2_{DB}(X,\mathbb{Z}(1))\to H^3(\overline{X}^{an},X^{an};\mathbb{Z})\to H^2(\overline{X}^{an},{\cal O}^*_{\overline{X}^{an}})$$
c) There is an exact sequence
$$0\to H^1(\overline{X}^{an},j_*^m{\cal O}^*_{X^{an}})\to H^2_{DB}(X,\mathbb{Z}(1))\to H^0(\overline{X}^{an},R^2j^{an}_*\mathbb{Z}_{X^{an}})$$

{\bf Proof:} a) As in Theorem 2.4 we have an exact sequence
\footnotesize
$$\begin{array}{ccccccc}
H^2(\overline{X}^{an},X^{an};\mathbb{Z})&\to& H^2_D(\overline{X},\mathbb{Z}(1))&\to&Pic\,X&\to& 0\\
\downarrow\simeq&&\downarrow\simeq&&\downarrow&&\downarrow\\
H^0(\overline{X}^{an},j_*^m{\cal O}^*_{X^{an}}/{\cal O}^*_{\overline{X}^{an}})&\to&Pic\,\overline{X}&\to& H^1(\overline{X}^{an},j_*^m{\cal O}^*_{X^{an}})&\to& H^1(\overline{X}^{an},j_*^m{\cal O}^*_{X^{an}}/{\cal O}^*_{\overline{X}^{an}})\\
\end{array}$$
\normalsize
This implies that $0\to Pic\,X\to H^1(\overline{X}^{an},j_*^m{\cal O}^*_{X^{an}})\to H^1(\overline{X}^{an},j_*^m{\cal O}^*_{X^{an}}/{\cal O}^*_{\overline{X}^{an}})$ is exact. The rest is clear because $j_*^m{\cal O}^*_{X^{an}}/{\cal O}^*_{\overline{X}^{an}}\simeq R^1j^{an}_*\mathbb{Z}_{X^{an}}\simeq R^2\Gamma_{D^{an}}\mathbb{Z}_{\overline{X}^{an}}\simeq \oplus^r_{j=1}(i_j)^{an}_*\mathbb{Z}_{D_j^{an}}$, where $i_j:D_j\to \overline{X}$ is the inclusion, see (*) of the proof of Theorem 2.4. Note that the lower exact sequence can be extended to the right by $\to H^2(\overline{X}^{an},{\cal O}^*_{\overline{X}^{an}})$.

\vskip.1in
b) The beginning of the proof of Theorem 2.4 shows the exactness of the sequence\\
$H^2_{DB}(X,\mathbb{Z}(1))\to H^3(\overline{X}^{an},X^{an};\mathbb{Z})\to H^2(\overline{X}^{an},{\cal O}^*_{\overline{X}^{an}})$.\\
The rest follows from Theorem 2.4 and Corollary 1.6.

\vskip.1in
c) Look at the following commutative diagram:
\footnotesize
$$\begin{array}{ccccccccc}
H^0(\overline{X},j_*^m{\cal O}^*_{X^{an}}/{\cal O}^*_{\overline{X}^{an}})&\to&Pic\,\overline{X}&\to&H^1(\overline{X},j_*^m{\cal O}^*_{X^{an}})&\to&H^1(\overline{X}^{an},R^2\Gamma_{D^{an}}\mathbb{Z}_{\overline{X}^{an}})\\
\downarrow\simeq&&\downarrow\simeq&&\downarrow&&\downarrow\\
H^2(\overline{X}^{an},X^{an};\mathbb{Z})&\to&Pic\,\overline{X}&\to&H^2_{DB}(\overline{X},X,\mathbb{Z}(1))&\to&H^3(\overline{X}^{an},X^{an};\mathbb{Z})
\end{array}$$
\normalsize
On the right we may extend the rows by $\to H^2(\overline{X},{\cal O}^*_{\overline{X}^{an}})$. Note that $j_*^m{\cal O}^*_{X^{an}}/{\cal O}^*_{\overline{X}^{an}}\simeq R^2\Gamma_{D^{an}}\mathbb{Z}_{\overline{X}^{an}}$, see above.\\
Since $R^k\Gamma_{D^{an}}\mathbb{Z}_{\overline{X}^{an}}=0, k=0,1$, we have an exact sequence\\
$0\to H^1(\overline{X},R^2\Gamma_{D^{an}}\mathbb{Z}_{\overline{X}^{an}})\to H^3(\overline{X}^{an},X^{an};\mathbb{Z})\to H^0(\overline{X},R^3\Gamma_{D^{an}}
\mathbb{Z}_{\overline{X}^{an}})$. This implies our statement; note that $R^3\Gamma_{D^{an}}\mathbb{Z}_{\overline{X}^{an}}\simeq R^2j^{an}_*\mathbb{Z}_{X^{an}}$.\\
Assume e.g. that $x\in H^2_{DB}(X,\mathbb{Z}(1))$ is mapped to $0\in H^0(\overline{X}^{an},R^2j^{an}_*\mathbb{Z}_{X^{an}})$. Then the image of $x$ in $H^3(\overline{X}^{an},X^{an};\mathbb{Z})$ comes from $y\in H^1(\overline{X}^{an},R^2\Gamma_{D^{an}}\mathbb{Z}_{\overline{X}^{an}})$. Since $x$ is mapped to $0\in H^2(\overline{X}^{an},{\cal O}^*_{\overline{X}^{an}})$, see b), we have that $y\mapsto 0\in H^2(\overline{X}^{an},{\cal O}^*_{\overline{X}^{an}})$, so $y$ comes from $z\in H^1(\overline{X}^{an},j^m_*{\cal O}^*_{X^{an}})$, by the upper line of the commutative diagram above. Then $z\mapsto x$.\\

{\bf Corollary 2.7:} We have injective mappings $Pic\,X\to H^1(\overline{X}^{an},j_*^m{\cal O}^*_{X^{an}})\to H^2_{DB}(X,\mathbb{Z}(1))$.\\

In general $H^1(\overline{X},j_*{\cal O}^*_X)\simeq Pic\,X\not\simeq H^1(\overline{X}^{an},j_*^m{\cal O}^*_{X^{an}})$
as Example b) in section 4 shows.\\

In particular we obtain the following description of $H^1(\overline{X}^{an},j_*^m{\cal O}^*_{X^{an}})$ from Proposition 2.6c):\\

{\bf Lemma 2.8:} $\{[L]\in H^2_{DB}(X,\mathbb{Z}(1))\,|\,c_{1,x}(L)=0\,\forall x\in D^{an}\}\simeq H^1(\overline{X},j_*^m{\cal O}^*_{X^{an}})$.\\
In particular, $H^2_{DB}(X,\mathbb{Z}(1))\simeq H^1(\overline{X}^{an},j_*^m{\cal O}^*_{X^{an}})$ if $D$ is smooth.\\

Note that an element of $H^2_{DB}(X,\mathbb{Z}(1))$ yields an element of $H^2_D(X^{an},\mathbb{Z}(1))\simeq Pic\,X^{an}$ which can be represented by a line bundle $L$. And $c_{1,x}(L)$ is the first Chern class of $L|U\setminus D^{an}$, $U$ being a convenient neighbourhood of $x$ in $\overline{X}^{an}$.\\

But in general the latter is no longer true, see Example a) in section 4.\\

We have a geometric interpretation of $H^1(\overline{X}^{an},j_*^m{\cal O}^*_{X^{an}})$:\\
According to [D] II Prop. 2.22, p. 70, we know that $Pic\,X$ coincides with the group of isomorphism classes of line bundles on $X^{an}$ which are effectively meromorphic along $D^{an}$. For the definition of meromorphic and effectively meromorphic line bundles see [D] II 2.13, p. 65f., or [M] p. 153f. Now\\

{\bf Remark 2.9:} $H^1(\overline{X}^{an},j_*^m{\cal O}^*_{X^{an}})$ coincides with the group of isomorphism classes of line bundles on $X^{an}$ which are meromorphic along $D^{an}$.\\

This will be proved in a subsequent paper.\\

{\bf Remark 2.10:} Since the injection $Pic\,X\to H^1(\overline{X}^{an},j_*^m{\cal O}^*_{X^{an}})$ is not bijective, in general, see above, we conclude that meromorphic line bundles are not necessarily effectively meromorphic. 

\vskip.1in
This answers a question by Deligne [D], p. 66, and Malgrange [M], p. 154.\\

{\bf Remark 2.11:} In the case of the algebraic Picard group there is no obvious reason why to restrict to the smooth case. In the singular case one can proceed as follows:\\

1. Suppose that $X$ is a complex algebraic variety, $Sing\,X$ compact. Define $H^k_{DB}$ and $H^2_{db}$ as before.\\
Then $Pic\, X\simeq H^2_{db}(X,\mathbb{Z}(1))$:\\
Proceed as in the proof of Theorem 2.4, using a compactification $\overline{X}$ of $X$ such that $\overline{X}\setminus Sing\,X$ is smooth and $D:=\overline{X}\setminus X$ is a divisor with normal crossings.\\
We have again that  $H^2_{db}(X,\mathbb{Z}(1))$ is the set of all elements of $H^2_{DB}(X,\mathbb{Z}(1))$ whose image in $H^2(X^{an};\mathbb{C})$ is contained in $W_2H^2(X^{an};\mathbb{C})$.\\
In fact, note that Lemma 0.3 still holds in our context: For the cohomology of $(\overline{X}^{an},X^{an})$ it does not matter whether $X$ is smooth or not, by excision.\\
Furthermore $W_2H^2(X^{an};\mathbb{C})=I$, where $I$ is the image of $H^2(\overline{X}^{an};\mathbb{C})\to H^2(X^{an};\mathbb{C})$:\\
First, $\overline{X}$ is compact, so $H^2(\overline{X}^{an};\mathbb{C})=W_2H^2(\overline{X}^{an};\mathbb{C})$, hence $I\subset W_2H^2(X^{an};\mathbb{C})$.\\
Now we have an exact sequence $H^2(\overline{X}^{an};\mathbb{C})\to H^2(X^{an};\mathbb{C})\to H^3_D(\overline{X}^{an};\mathbb{C})$, as in [E] Prop. 3.7.18, and the weights of the mixed Hodge structure on $H^3_D(\overline{X}^{an};\mathbb{C})$ vary between $3$ and $4$, see loc.cit. (by excision assuming $X$ to be smooth). So $W_2H^3_D(\overline{X}^{an};\mathbb{C})=0$.
Since morphisms of mixed Hodge structures are strict we obtain an exact sequence $H^2(\overline{X}^{an};\mathbb{C})\to W_2H^2(X^{an};\mathbb{C})\to 0$, so $W_2H^2(X^{an};\mathbb{C})\subset I$.\\

2. If $X$ is smooth, $Pic\,X$ coincides with the group of classes of Cartier divisors as well as of Weil divisors. This holds no longer in general.\\
Now suppose that ${\rm codim}\,Sing\,X\ge 2$ (e.g. $X$ normal). Let $Cl\,X$ be the group of classes of Weil divisors on $X$. Then $Cl\,X\simeq H^2_{db}(X\setminus Sing\,X,\mathbb{Z}(1))$:\\
According to [Ha] II Prop. 6.5, p. 133, we have $Cl\,X\simeq Cl\,X\setminus Sing\,X$, the rest follows from Theorem 2.4.\\

{\bf Proposition 2.12:} Let $f:X\to Y$ be a morphism between smooth complex algebraic varieties. Assume that $f$ induces isomorphisms $H^k(Y^{an};\mathbb{Z})\simeq H^k(X^{an};\mathbb{Z})$, $k=1,2$. Then:\\
a) $H^2_{DB}(Y,\mathbb{Z}(1))\simeq H^2_{DB}(X,\mathbb{Z}(1))$,\\
b) (see [HL1] Theorem 1.4): $Pic\,Y\simeq Pic\,X$.\\

{\bf Proof:} a) This follows from the commuative diagram with exact rows:
\footnotesize
$$\begin{array}{ccccccccc}
H^1(Y^{an};\mathbb{Z})&\to&H^1(\overline{Y}^{an},{\cal O}_{\overline{Y}^{an}})&\to&H^2_{DB}(Y,\mathbb{Z}(1))&\to&H^2(Y^{an};\mathbb{Z})&\to&H^2(\overline{Y}^{an},{\cal O}_{\overline{Y}^{an}})\\
\downarrow&&\downarrow&&\downarrow&&\downarrow&&\downarrow\\
H^1(X^{an};\mathbb{Z})&\to&H^1(\overline{X}^{an},{\cal O}_{\overline{X}^{an}})&\to&H^2_{DB}(X,\mathbb{Z}(1))&\to&H^2(X^{an};\mathbb{Z})&\to&H^2(\overline{X}^{an},{\cal O}_{\overline{X}^{an}})
\end{array}$$
\normalsize
Here we take convenient compactifications of $X, Y$ such that $f$ extends to a morphism $\overline{X}\to\overline{Y}$.\\
Note that $H^k(\overline{X}^{an},{\cal O}_{\overline{X}^{an}})\simeq Gr^0_FH^k(X^{an};\mathbb{C})$. Therefore the first, second,fourth and fifth vertical are isomorphisms, so the third one, too.\\
Or use Lemma 2.15 below.

\vspace{2mm}
b) This follows from a) because $Pic\,X\simeq H^2_{db}(X,\mathbb{Z}(1))$.\\
In fact, the injectivity is clear because $H^2_{db}\subset H^2_{DB}$ and $H^2_{DB}(Y,\mathbb{Z}(1))\to H^2_{DB}(X,\mathbb{Z}(1))$ is injective.\\
Surjectivity: Assume that $c\in H^2_{db}(X,\mathbb{Z}(1))$, i.e. $c\in H^2_{DB}(X,\mathbb{Z}(1))$ and the image in $H^2(X^{an};\mathbb{C})$ is contained in $W_2H^2(X^{an};\mathbb{C})$. Then $c$ has an inverse image in $H^2_{DB}(Y,\mathbb{Z}(1))$, the image in $H^2(Y^{an};\mathbb{C})$ is mapped onto the image of $c$, so it is contained in $W_2H^2(Y^{an};\mathbb{C})$. 

\vspace{2mm}
Alternatives: use Lemma 2.14 or 2.16 below.\\

We can define an analogue of the N\'eron-Severi group $NS(X):=Im(Pic\,X\to H^2(X^{an};\mathbb{Z}))$ for $H^2_{DB}(X,\mathbb{Z}(1))$ instead of $Pic\,X$:\\
$NS_{DB}(X):=Im(H^2_{DB}(X,\mathbb{Z}(1))\to H^2(X^{an};\mathbb{Z}))$.\\

{\bf Lemma 2.13:} $NS_{DB}(X)=\gamma^{-1}(F^1H^2(X^{an};\mathbb{C}))$, where $\gamma:H^2(X^{an};\mathbb{Z})\to H^2(X^{an};\mathbb{C})$.\\

{\bf Proof:} The exact sequence $H^2_{DB}(X,\mathbb{Z}(1))\to H^2(X^{an};\mathbb{Z})\to H^2(\overline{X}^{an},{\cal O}_{\overline{X}^{an}})$ shows that $NS_{DB}(X)=ker(H^2(X^{an};\mathbb{Z})\to H^2(\overline{X}^{an},{\cal O}_{\overline{X}^{an}})$. Now note that $F^1H^2(X^{an};\mathbb{C})=ker(H^2(X^{an};\mathbb{C})\to H^2(\overline{X}^{an},{\cal O}_{\overline{X}^{an}}))$.\\

Compare with [HL1] Theorem 3.1: $NS(X)=\gamma^{-1}(F^1W_2H^2(X^{an};\mathbb{C}))$.\\

We have a precise description of the difference between $Pic\,X$ and $H^2_{DB}(X,\mathbb{Z}(1))$:\\

{\bf Lemma 2.14:} We have an exact sequence
$$0\to Pic\,X\to H^2_{DB}(X,\mathbb{Z}(1))$$
$$\to (H^2(X^{an};\mathbb{Z})/W_2H^2(X^{an},\mathbb{Z}))\cap (F^1H^2(X^{an};\mathbb{C})/F^1W_2H^2(X^{an};\mathbb{C}))\to 0$$
Note that we may consider $H^2(X^{an};\mathbb{Z})/W_2H^2(X^{an},\mathbb{Z})$ as a subgroup of $H^2(X^{an};\mathbb{C})/W_2H^2(X^{an},\mathbb{C})$ because the natural mapping $H^2(X^{an};\mathbb{Z})/W_2H^2(X^{an},\mathbb{Z})\to H^2(X^{an};\mathbb{C})/W_2H^2(X^{an},\mathbb{C})$ is injective.\\

{\bf Proof:} It is sufficient to show that
$$0\to Pic\,X/Pic^0\,X\to H^2_{DB}(X,\mathbb{Z}(1))/Pic^0\,X$$
$$\to (H^2(X^{an};\mathbb{Z})/W_2H^2(X^{an},\mathbb{Z}))\cap (F^1H^2(X^{an};\mathbb{C})/F^1W_2H^2(X^{an};\mathbb{C}))\to 0$$
is exact, i.e. that the sequence
$$0\to \gamma^{-1}(F^1W_2H^2(X^{an};\mathbb{C}))\to \gamma^{-1}(F^1H^2(X^{an};\mathbb{C}))$$
$$\to (H^2(X^{an};\mathbb{Z})/W_2H^2(X^{an},\mathbb{Z}))\cap (F^1H^2(X^{an};\mathbb{C})/F^1W_2H^2(X^{an};\mathbb{C}))\to 0$$
is exact, using Lemma 2.13 and [HL1] loc. cit. But this is easy to verify.\\

Replacing $\mathbb{Z}$ by $\mathbb{Q}$, we may deduce that the cokernel of the mapping $H^2_{AH}(X;\mathbb{Z}(1))\to H^2_{DB}(X;\mathbb{Z}(1))$ is not finite if and only if $(H^2(X^{an};\mathbb{Q})/W_2H^2(X^{an},\mathbb{Q}))\cap (F^1H^2(X^{an};\mathbb{C})/F^1W_2H^2(X^{an};\mathbb{C}))\neq 0$, cf. [S] Remark (3.5) (i), p. 297.\\

We have an exact sequence
$$H^1(X^{an};\mathbb{Z})\to H^1(\overline{X}^{an},{\cal O}_{\overline{X}^{an}})\to H^2_{DB}(X,\mathbb{Z}(1))\to H^2(X^{an};\mathbb{Z})\to H^2(\overline{X}^{an},{\cal O}_{\overline{X}^{an}})$$
This leads to a short exact sequence:\\

{\bf Lemma 2.15:} There is an exact sequence
$$0\to coker(H^1(X^{an};\mathbb{Z})\to H^1(X^{an};\mathbb{C})/F^1H^1(X^{an};\mathbb{C}))\to H^2_{DB}(X,\mathbb{Z}(1))$$
$$\to ker(H^2(X^{an};\mathbb{Z})\to H^2(X^{an};\mathbb{C})/F^1H^2(X^{an};\mathbb{C}))\to 0$$
(see [S] (3.2.1), p. 292).\\

{\bf Proof:} We have a long exact sequence
$$\ldots\to H^k(\overline{X},\Omega^{\ge 1}_{\overline{X}^{an}}(log\,D))\to H^k(\overline{X},\Omega^\cdot_{\overline{X}^{an}}(log\,D))\to H^k(\overline{X},{\cal O}_{\overline{X}^{an}})\to\ldots$$
The first mapping corresponds to the inclusion $F^1H^k(X^{an};\mathbb{C})\to H^k(X^{an};\mathbb{C})$, so it is injective. Therefore we obtain short exact sequences, and $H^k(\overline{X},{\cal O}_{\overline{X}^{an}})\simeq H^k(X^{an};\mathbb{C})/F^1H^k(X^{an};\mathbb{C})$. The rest is clear from the preceding long exact sequence.\\

Now we pass to absolute Hodge cohomology:\\

{\bf Lemma 2.16:} There is an exact sequence
\footnotesize
$$0\to coker(W_2H^1(X^{an};\mathbb{Z})\to W_2H^1(X^{an};\mathbb{C})/F^1W_2H^1(X^{an};\mathbb{C}))\to H^2_{AH}(X,\mathbb{Z}(1))$$
$$\to ker(W_2H^2(X^{an};\mathbb{Z})\to W_2H^2(X^{an};\mathbb{C})/F^1W_2H^2(X^{an};\mathbb{C}))\to 0$$
\normalsize
Of course we may replace $H^2_{AH}(X,\mathbb{Z}(1))$ by $Pic\,X$ here, see Theorem 0.4.\\

{\bf Proof:} According to [S] (3.2.2), p. 292, we have an exact sequence
\footnotesize
$$0\to Ext^1_{MHS(\mathbb{Z})}(\mathbb{Z},H^1(X^{an};\mathbb{Z}(1))\to H^2_{AH}(X,\mathbb{Z}(1))$$
$$\to Hom_{MHS(\mathbb{Z})}(\mathbb{Z},H^2(X^{an};\mathbb{Z}(1))\to 0$$
\normalsize
If we take $\mathbb{Q}$ instead of $\mathbb{Z}$ this implies our lemma with $\mathbb{Q}$ instead of $\mathbb{Z}$, see [S] (3.2.4), p. 293. But in the case of the ring $\mathbb{Z}$ we proceed directly:\\
First note that $W_2H^1(X^{an};\mathbb{Z})=H^1(X^{an};\mathbb{Z})$, so $coker(W_2H^1(X^{an};\mathbb{Z})\to W_2H^1(X^{an};\mathbb{C})/F^1W_2H^1(X^{an};\mathbb{C}))$ may be simpler written as $coker(H^1(X^{an};\mathbb{Z})\to H^1(X^{an};\mathbb{C})/F^1H^1(X^{an};\mathbb{C}))$.\\
Now we look at the exact sequence 
$$0\to Pic^0\,X\to Pic\,X\to NS\,X\to 0$$
First, $Pic^0\,X\simeq ker(H^2_{DB}(X,\mathbb{Z}(1))\to H^2(X^{an};\mathbb{Z}))\simeq coker(H^1(X^{an};\mathbb{Z})\to H^1(X^{an};\mathbb{C})/F^1H^1(X^{an};\mathbb{C}))$, cf. Lemma 2.15, and $NS\,X\simeq ker(W_2H^2(X^{an};\mathbb{Z})\to W_2H^2(X^{an};\mathbb{C})/F^1W_2H^2(X^{an};\mathbb{C}))$, by [HL1] Theorem 3.1.\\

{\bf 3. Algebraic line bundles with connections}\\

Let $Pic_{cir}\,X$ be the group of isomorphism classes of line bundles on $X$ with a regular integrable connection.\\

{\bf Theorem 3.1:} a) $H^2_{db}(X,\mathbb{Z}(2))\simeq Pic_{cir}\,X$.\\
b) $H^2_{DB}(X,\mathbb{Z}(p))\simeq H^2_{db}(X,\mathbb{Z}(p))\simeq Pic_{cir}\,X, p\ge 3$.\\

{\bf Proof:} a) Look at $p=2$. Then we have an exact sequence:
\footnotesize
$$H^1_{DB}(X,\mathbb{Z}(1))\to H^0(\overline{X}^{an},\Omega^1_{\overline{X}^{an}}(log\,D))\to H^2_{DB}(X,\mathbb{Z}(2))\to H^2_{DB}(X,\mathbb{Z}(1))\to H^1(\overline{X}^{an},\Omega^1_{\overline{X}^{an}}(log\,D))$$
\normalsize
which induces an exact sequence
\footnotesize
$$H^1_{DB}(X,\mathbb{Z}(1))\to H^0(\overline{X}^{an},\Omega^1_{\overline{X}^{an}}(log\,D))\to H^2_{db}(X,\mathbb{Z}(2))\to H^2_{db}(X,\mathbb{Z}(1))\to H^1(\overline{X}^{an},\Omega^1_{\overline{X}^{an}}(log\,D))$$
\normalsize
In fact, the second mapping is well-defined because $d:H^2_{DB}(X,\mathbb{Z}(2))\to H^3(\overline{X}^{an},X^{an};\mathbb{Z})$ factorizes over $H^2_{DB}(X,\mathbb{Z}(1))$.\\
On the other hand, we have an exact sequence
$$H^0(X,{\cal O}^*_X)\to H^0(\overline{X},\Omega^1_{\overline{X}}(log\,D))\to Pic_{cir}(X)\to Pic(X)\to H^1(\overline{X},\Omega^1_{\overline{X}}(log\,D))$$ 
see [HL2] Theorem 3.10, Lemma 3.11.\\
So $Pic_{cir}(X)\simeq H^2_{db}(X,\mathbb{Z}(2))$, as soon as the two exact sequences fit together to a commutative diagram
\footnotesize
$$\begin{array}{cccccc}
H^0(X,{\cal O}^*_X)&\to& H^0(\overline{X},\Omega^1_{\overline{X}}(log\,D))&\to& Pic_{cir}(X)&\to\\
\downarrow\simeq&&\downarrow\simeq&&\downarrow&\\
H^1_{DB}(X,\mathbb{Z}(1))&\to& H^0(\overline{X}^{an},\Omega^1_{\overline{X}^{an}}(log\,D))&\to& H^2_{db}(X,\mathbb{Z}(2))&\to
\end{array}$$
$$\begin{array}{cccc}
\to& Pic(X)&\to& H^1(\overline{X},\Omega^1_{\overline{X}}(log\,D))\\
&\downarrow\simeq&&\downarrow\simeq\\
\to& H^2_{db}(X,\mathbb{Z}(1))&\to& H^1(\overline{X}^{an},\Omega^1_{\overline{X}^{an}}(log\,D))
\end{array}$$
\normalsize
In order to show this, let us start with the commutative diagram
\footnotesize
$$\begin{array}{cccccc}
H^0(X,{\cal O}^*_X)&\to& H^0(\overline{X},\Omega^1_{\overline{X}}(log\,D))&\to& Pic_{cir}(X)&\to\\
\downarrow\simeq&&\downarrow\simeq&&\downarrow&\\
H^0(\overline{X}^{an},j_*^m{\cal O}^*_{X^{an}})&\to& H^0(\overline{X}^{an},\Omega^1_{\overline{X}^{an}}(log\,D))&\to& \mathbb{H}^1({\cal T}^\cdot)&\to
\end{array}$$
$$\begin{array}{cccc}
\to& Pic(X)&\to& H^1(\overline{X},\Omega^1_{\overline{X}}(log\,D))\\
&\downarrow&&\downarrow\simeq\\
\to&H^1(\overline{X}^{an},j_*^m{\cal O}^*_{X^{an}})&\to&H^1(\overline{X}^{an},\Omega^1_{\overline{X}^{an}}(log\,D))
\end{array}$$
\normalsize
where ${\cal T}^\cdot:=cone(j_*^m{\cal O}^*_{X^{an}}\to \Omega^1_{\overline{X}^{an}}(log\,D))[-1]$. Note that $Pic_{cir}\simeq \mathbb{H}^1(\overline{X},{\cal S}^\cdot)$, where ${\cal S}^\cdot:=cone(j_*{\cal O}^*_X\to \Omega^1_X(log\,D))[-1]$. See [HL2] \S 3.\\
The first, second, and fifth vertical are bijective by GAGA.\\
The lower exact sequence is the long exact cohomology sequence for
$$0\to \Omega^1_{\overline{X}^{an}}(log\,D)[1]\to cone (\tilde{\mathbb{Z}}(1)[-1]\to \Omega^1_{\overline{X}^{an}}(log\,D))\to \tilde{\mathbb{Z}}(1)\to 0$$
cf. Lemma 2.3. This sequence can be mapped to
$$0\to \Omega^1_{\overline{X}^{an}}(log\,D)[1]\to cone (\mathbb{Z}(1)_{DB}[-1]\to \Omega^1_{\overline{X}^{an}}(log\,D))\to \mathbb{Z}(1)_{DB}\to 0$$
Now the cone in the middle is quasiisomorphic to $\mathbb{Z}(2)_{DB}$.\\
So we obtain a commutative diagram
\footnotesize
$$\begin{array}{cccccc}
H^0(X,{\cal O}^*_X)&\to& H^0(\overline{X},\Omega^1_{\overline{X}}(log\,D))&\to& Pic_{cir}(X)&\to\\
\downarrow\simeq&&\downarrow\simeq&&\downarrow&\\
H^1_{DB}(X,\mathbb{Z}(1))&\to&H^0(\overline{X}^{an},\Omega^1_{\overline{X}^{an}}(log\,D))&\to& H^2_{DB}(X,\mathbb{Z}(2))&\to
\end{array}$$
$$\begin{array}{cccc}
\to&Pic(X)&\to& H^1(\overline{X},\Omega^1_{\overline{X}}(log\,D))\\
&\downarrow&&\downarrow\simeq\\
\to&H^2_{DB}(X,\mathbb{Z}(1))&\to&H^1(\overline{X}^{an},\Omega^1_{\overline{X}^{an}}(log\,D))
\end{array}$$
\normalsize
The left vertical is bijective, see Lemma 0.3 (i.e. [EV]), the second and fifth vertical are bijective, too, as we know already.\\
Now we may replace $H^2_{DB}$ by $H^2_{db}$, see above. Then $Pic\;X\to H^2_{db}(X,\mathbb{Z}(1))$ is bijective by Theorem 2.4. The rest is clear.\\

b) Similarly, we have an exact sequence
$$0\to H^2_{DB}(X,\mathbb{Z}(3))\to H^2_{DB}(X,\mathbb{Z}(2))\to H^0(\overline{X},\Omega^2_{\overline{X}}(log\,D))$$
which induces an exact sequence
$$0\to H^2_{db}(X,\mathbb{Z}(3))\to H^2_{db}(X,\mathbb{Z}(2))\to H^0(\overline{X},\Omega^2_{\overline{X}}(log\,D))$$
But the elements of $H^2_{db}(X,\mathbb{Z}(2))$ correspond to certain connections which are integrable, hence $H^2_{db}(X,\mathbb{Z}(3))\simeq H^2_{db}(X,\mathbb{Z}(2))$, because the last arrow is given by the curvature of the connection in question.\\
Now for $p\ge 3$, $H^2_{DB}(X,\mathbb{Z}(p+1))\simeq H^2_{DB}(X,\mathbb{Z}(p))$, hence $H^2_{db}(X,\mathbb{Z}(p+1))\simeq H^2_{db}(X,\mathbb{Z}(p))$.\\
Altogether, $H^2_{db}(X,\mathbb{Z}(p))\simeq H^2_{db}(X,\mathbb{Z}(2))\simeq Pic_{cir}(X)$, $p\ge 2$.\\
Finally $H^2_{DB}(X,\mathbb{Z}(p))=H^2_{db}(X,\mathbb{Z}(p))$ for $p\ge 3$: This can be seen in different ways.

\vskip.1in
(i) Let $a\in H^2_{DB}(X,\mathbb{Z}(p))$. The image of $a$ in $H_D^2(X^{an},\mathbb{Z}(p))\simeq Pic_{ci}(X^{an})$ is represented by an analytic line bundle with an integrable connection, so its first complex Chern class vanishes. But this is the image of $c(a)$. By Lemma 1.5 we obtain our statement.\\
(ii) It is sufficient to show this for $p\ge \dim\;X$. Then $H^2_{DB}(X,\mathbb{Z}(p))$ is the hypercohomology of the cone of $Rj^{an}_*\mathbb{Z}_{X^{an}}\to \Omega^\cdot_{\overline{X}^{an}}(log\,D)$. Now $\Omega^\cdot_{\overline{X}^{an}}(log\,D)$ is quasiisomorphic to $j^{an}_*\Omega^\cdot_{X^{an}}$, hence $H^2_{DB}(X,\mathbb{Z}(p))\simeq H^2_D(X^{an},\mathbb{Z}(p))\simeq Pic_{ci}(X^{an})\simeq Pic_{cir}(X)\simeq H^2_{db}(X,\mathbb{Z}(p))$ by Deligne's existence theorem, see [D] II Th. 5.9, p. 97.\\
(iii) As in (i), we have a natural map $H^2_{DB}(X,\mathbb{Z}(p))\to  H^2_D(X^{an},\mathbb{Z}(p))\simeq Pic_{ci}^{an}(X)$.
Now the natural map $H^2_D(X^{an},\mathbb{Z}(p))\to H^3(\overline{X}^{an},X^{an};\mathbb{Z})$ is the zero map, because of Deligne's existence theorem which guarantees that we deal with algebraic line bundles.\\

In order to take connections into account which are not integrable or which are irregular, we modify the definition of Deligne-Beilinson cohomology: Let ${\cal S}(p):=cone(Rj_*^{an}\mathbb{Z}(p)\to ({\cal O}_{\overline{X}^{an}}\to j_*^m\Omega^1_{X^{an}}\to\ldots\to j_*^m\Omega^{p-1}_{X^{an}}\to 0\to \ldots))[-1]$. Let 
$$H^k_{DH}(X,\mathbb{Z}(p)):=\mathbb{H}^k(\overline{X}^{an},{\cal S}(p)) \hbox{  (hybrid Deligne cohomology) and} $$ 
$$H^2_{dh}(X,\mathbb{Z}(p)):=ker(d^*:H^2_{DH}(X,\mathbb{Z}(p))\to H^3(\overline{X}^{an},X^{an};\mathbb{Z})$$ 
Here $d^*$ is defined as the composition $H^2_{DH}(X,\mathbb{Z}(p))\to H^2_{DH}(X,\mathbb{Z}(0)))\simeq H^2(X^{an};\mathbb{Z})\to H^3(\overline{X}^{an},X^{an};\mathbb{Z})$. \\
Since we know that $H^2_{DB}(X,\mathbb{Z}(1))$ is independent of the compactification we see that $H^2_{DH}(X,\mathbb{Z}(p))$ and $H^2_{dh}(X,\mathbb{Z}(p))$ are independent, too.\\
Note that $H^2_{dh}(X,\mathbb{Z}(1))\simeq H^2_{db}(X,\mathbb{Z}(1))\simeq Pic\,X$.\\

{\bf Theorem 3.2:} a) $H^2_{dh}(X,\mathbb{Z}(2))\simeq Pic_c\,X$.\\
b) $H^2_{DH}(X,\mathbb{Z}(p))\simeq H^2_{dh}(X,\mathbb{Z}(p))\simeq Pic_{ci}\,X, p\ge 3$.\\

{\bf Proof:} a) There is an exact sequence
$$H^1_{DH}(X,\mathbb{Z}(1))\to H^0(\overline{X}^{an},j_*^m\Omega^1_{X^{an}})\to H^2_{DH}(X,\mathbb{Z}(2))\to H^2_{DH}(X,\mathbb{Z}(1))\to H^1(\overline{X}^{an},j_*^m\Omega^1_{X^{an}})$$
which leads to an exact sequence
$$H^1_{DH}(X,\mathbb{Z}(1))\to H^0(\overline{X}^{an},j_*^m\Omega^1_{X^{an}})\to H^2_{dh}(X,\mathbb{Z}(2))\to H^2_{dh}(X,\mathbb{Z}(1))\to H^1(\overline{X}^{an},j_*^m\Omega^1_{X^{an}})$$
Again the second mapping is well-defined because $d^*:H^2_{DH}(X,\mathbb{Z}(2))\to H^3(\overline{X}^{an},X^{an};\mathbb{Z})$ factors over $H^2_{DH}(X,\mathbb{Z}(1))$.\\
On the other hand, there is an exact sequence
$$H^0(X,{\cal O}^*_X)\to H^0(X,\Omega^1_X)\to Pic_c(X)\to Pic(X)\to H^1(X,\Omega^1_X)$$\\
see [HL2] Theorem 3.1, so $Pic_c(X)\simeq H^2_{dh}(X,\mathbb{Z}(2))$. Note that $H^1_{DH}(X,\mathbb{Z}(1))= H^1_{DB}(X,\mathbb{Z}(1))\simeq H^0(X,{\cal O}^*_X)$, cf. Lemma 0.3.\\
In detail we proceed similarly as in the case of Theorem 3.1, with $j_*\Omega^1_X$ instead of $\Omega^1_X(log\,D)$. Note that $H^p(X,\Omega^1_X)\simeq H^p(\overline{X}^{an},j_*^m\Omega^1_{X^{an}})$ by GAGA.\\

b) Finally, the exact sequence
$$0\to H^2_{DH}(X,\mathbb{Z}(3))\to H^2_{DH}(X,\mathbb{Z}(2))\to H^0(\overline{X}^{an},j_*^m\Omega^2_{X^{an}})$$
leads to an exact sequence
$$0\to H^2_{dh}(X,\mathbb{Z}(3))\to H^2_{dh}(X,\mathbb{Z}(2))\to H^0(\overline{X}^{an},j_*^m\Omega^2_{X^{an}})$$
Comparison with the exact sequence
$$0\to Pic_{ci}(X)\to Pic_c(X)\to H^0(X,\Omega^2_X)$$
gives that $Pic_{ci}(X)\simeq H^2_{dh}(X,\mathbb{Z}(3))$.\\
More exactly: there is a commutative diagram
$$\begin{array}{ccccccc}
0&\to&Pic_{ci}\,X&\to& Pic_c\,X&\to&H^0(X,\Omega^2_X)\\
&&\downarrow&&\downarrow\simeq&&\downarrow\simeq\\
0&\to&H^2_{dh}(X,\mathbb{Z}(3))&\to& H^2_{dh}(X,\mathbb{Z}(2))&\to&H^0(\overline{X}^{an},j_*^m\Omega^2_{X^{an}})
\end{array}$$
hence $Pic_{ci}\,X\simeq H^2_{dh}(X,\mathbb{Z}(3))$.\\
In order to show this, recall first that $Pic_{ci}\,X=\mathbb{H}^1(X,{\cal S}^\cdot_0)$, where ${\cal S}^\cdot_0$ is the complex ${\cal O}^*_X\to {^c\Omega}^1_X\to 0\to\ldots$ with
${^c\Omega}^1_X:= ker(d:\Omega^1_X\to\Omega^2_X)$. See [HL2] Theorem 3.1.\\
Then ${\cal S}^\cdot_0$ is a subcomplex of ${\cal S}^\cdot$: ${\cal O}^*_X\to \Omega^1_X\to \Omega^2_X\to 0\to \ldots$, the quotient being $\overline{\cal S}^\cdot$: 
$0\to \Omega^1_X/{^c\Omega}^1_X\to \Omega^2_X\to 0\to\ldots$ which is quasiisomorphic to $0\to 0\to \Omega^2_X/d\Omega^1_X\to 0\to\ldots$.\\
So $\mathbb{H}^k(X,\overline{\cal S}^\cdot)=0$, $k=0,1$, hence $\mathbb{H}^1(X,{\cal S}_0^\cdot)\simeq \mathbb{H}^1(X,{\cal S}^\cdot)$.\\
Altogether, $0\to Pic_{ci}\,X\to Pic_c\,X\to H^0(X,\Omega^2_X)$ corresponds to the first line of the following commutative diagram:
$$\begin{array}{ccccccc}
0&\to&\mathbb{H}^1(\overline{X},j_*{\cal S}_3^\cdot)&\to&\mathbb{H}^1(\overline{X},j_*{\cal S}_2^\cdot)&\to&H^0(X,\Omega^2_X)\\
&&\downarrow&&\downarrow&&\downarrow\\
0&\to&\mathbb{H}^1(\overline{X}^{an},j_*^m({\cal T}^{an}_3)^\cdot)&\to&\mathbb{H}^1(\overline{X}^{an},j_*^m({\cal T}^{an}_2)^\cdot)&\to&H^0(\overline{X}^{an},j_*^m\Omega^2_{X^{an}})\\
&&\downarrow&&\downarrow&&\downarrow\\
0&\to&H^2_{dh}(X,\mathbb{Z}(3))&\to& H^2_{dh}(X,\mathbb{Z}(2))&\to&H^0(\overline{X}^{an},j_*^m\Omega^2_{X^{an}})
\end{array}$$
Here ${\cal S}_p^\cdot$ and ${\cal T}_p^\cdot$ are the complexes
${\cal O}^*_X\to\Omega^1_X\to\ldots\to\Omega^{p-1}_X\to 0\to\ldots$ resp.
${\cal O}^*_{X^{an}}\to\Omega^1_{X^{an}}\to\ldots\to\Omega^{p-1}_{X^{an}}\to 0\to\ldots$.\\
As for the transition from the second to the third line, note that, with $\Omega^{[1,p-1]}_{X^{an}}:=\Omega^1_{X^{an}}\to\ldots\to \Omega^{p-1}_{X^{an}}$:
$$j_*^m{\cal T}_p^{an}\sim cone(j_*^m{\cal O}_{X^{an}}^*\to j_*^m\Omega^{[1,p-1]}_{X^{an}})\sim cone(\tilde{\mathbb{Z}}(1)_X\to j_*^m\Omega^{[1,p-1]}_{X^{an}})$$
and we have a morphism 
$$cone(\tilde{\mathbb{Z}}(1)_X\to j_*^m\Omega^{[1,p-1]}_{X^{an}})\to cone(\mathbb{Z}(1)_X\to j_*^m\Omega^{[1,p-1]}_{X^{an}})$$
The composed mapping $\mathbb{H}^1(\overline{X},j_*{\cal S}_2^\cdot)\to H^2_{dh}(X,\mathbb{Z}(2))$ is bijective by a), and $H^0(X,\Omega^2_X)\to H^0(\overline{X}^{an},j_*^m\Omega^2_{X^{an}})$ is bijective by GAGA, hence $\mathbb{H}^1(\overline{X},j_*{\cal S}_3^\cdot)\to H^2_{dh}(X,\mathbb{Z}(3))$ is bijective, too.\\
It is easy to see that $H^2_{DH}(X,\mathbb{Z}(p+1))\simeq H^2_{DH}(X,\mathbb{Z}(p))$ for $p\ge 3$,\\
so $H^2_{dh}(X,\mathbb{Z}(p))\simeq H^2_{dh}(X,\mathbb{Z}(3))\simeq Pic_{ci}(X)$, $p\ge 3$.\\
Similarly, $H^2_{DH}(X,\mathbb{Z}(p))\simeq H^2_{dh}(X,\mathbb{Z}(p))$ for $p\ge 3$: We have a natural map $H^2_{DH}(X,\mathbb(p))\to  H^2_D(X^{an},\mathbb{Z}(p))\simeq Pic_{ci}(X^{an})$. Now we can argue similarly as in the proof of Theorem 3.1 b), using alternative (i) or (iii).\\

{\bf Consequence:} In comparison with the analytic case, the new objects are $Pic(X),$ $Pic_c(X)$ and $Pic_{ci}(X)$. Already $Pic(X)$ requires a modification of Deligne-Beilinson cohomology, this holds even more for $Pic_c(X)$ and $Pic_{ci}(X)$.\\ 

{\bf Remark 3.3:} One may ask whether it is more reasonable to look at $H^2_{D+}(X,\mathbb{Z}(p))=\mathbb{H}^2(\overline{X}^{an},cone(Rj_*^{an}\mathbb{Z}(p)\to j_*^m\Omega^{\le p}_{X^{an}})[-1])$ instead of $H^2_{DH}(X,\mathbb{Z}(p))$, and define $H^2_{d+}(X,\mathbb{Z}(p))$ similarly. But because of Grothendieck's comparison theorem [G], $H^2_{D+}(X,\mathbb{Z}(p))\simeq H^2_{\cal D}(X^{an},\mathbb{Z}(p))\simeq Pic_{ci}^{an}\,X^{an}\simeq Pic_{cir}\,X$ for $p=\dim\,X$, hence for $p\ge 3$.\\

Example b) of section 4 shows that $H^2_{D+}(X,\mathbb{Z}(1))\not\simeq H^2_{DB}(X,\mathbb{Z}(1))$, $H^2_{d+}(X,\mathbb{Z}(1))\not\simeq Pic\,X $, in general.\\

{\bf 4. Examples}\\

a) Put $X:=\mathbb{C}^*\times\mathbb{C}^*$. Let us compute the Deligne-Beilinson cohomology group $H^2_{DB}(X,\mathbb{Z}(1))$. Put $\overline{X}:=\mathbb{P}_1\times\mathbb{P}_1$. We have an exact sequence
$$H^1(\overline{X}^{an},{\cal O}_{\overline{X}^{an}})\to H^2_{DB}(X,\mathbb{Z}(1))\to H^2(X^{an},\mathbb{Z})\to H^2(\overline{X}^{an},{\cal O}_{\overline{X}^{an}})$$
But $H^k(\overline{X}^{an},{\mathcal O}_{\overline{X}^{an}})=0, k=1,2$. So the exact sequence is
$$0\to H^2_{DB}(X,\mathbb{Z}(1))\to \mathbb{Z}\to 0$$
i.e. $H^2_{DB}(X,\mathbb{Z}(1))\simeq \mathbb{Z}$.\\
On the other hand, $Pic\,X=0$ because $X\subset\mathbb{C}^2$ and $Pic\,\mathbb{C}^2=0$, see [Ha] II Prop. 6.2.\\
In particular, $d:H^2_{DB}(X,\mathbb{Z}(1))\to H^3(\overline{X}^{an},X^{an};\mathbb{Z})$ is not $0$.\\
Furthermore $H^0(X^{an},\Omega^2_{\overline{X}^{an}}(log\,D))\simeq \mathbb{C}$ with $\frac{dz_1}{z_1}\wedge \frac{dz_2}{z_2}$ as generator. Since $H^2(X^{an};\mathbb{C})\simeq\mathbb{C}$ we obtain by degeneration of the Hodge spectral sequence that $H^1(\overline{X}^{an},\Omega^1_{\overline{X}^{an}}(log\,D))=0$. Therefore $H^2_{DB}(X,\mathbb{Z}(2))\to H^2_{DB}(X,\mathbb{Z}(1))$ is surjective, so $d:H^2_{DB}(X,\mathbb{Z}(2))\to H^3(\overline{X}^{an},X^{an};\mathbb{Z})$ is not $0$.\\
Therefore $Pic_{cir}(X)\stackrel{\not\simeq}{\to}H^2_{DB}(X,\mathbb{Z}(2))$.\\
Similarly, $H^1(\overline{X}^{an},j_*^m\Omega^1_{X^{an}})=H^1(X,\Omega^1_X)=0$, so $H^2_{DH}(X,\mathbb{Z}(2))\to H^2_{DH}(X,\mathbb{Z}(1))$ is surjective, hence $d:H^2_{DH}(X,\mathbb{Z}(2))\to H^3(\overline{X}^{an},X^{an};\mathbb{Z})$ is not $0$, which implies $Pic_c(X)\stackrel{\not\simeq}{\to}H^2_{DH}(X,\mathbb{Z}(2))$.\\
So we cannot improve Theorem 2.4, 3.1a) and 3.2a).\\
Note that $0=Pic\,X\simeq H^1(\overline{X}^{an},j_*^m{\cal O}^*_{X^{an}})$: In our case $D=D_1\cup\ldots D_4$, and $D_j\simeq \mathbb{P}_1, j=1,\ldots,4$. Therefore $H^1(D_j^{an};\mathbb{Z})=0, j=1,\ldots,4$, so our assertion follows from Proposition 2.6.\\
So $H^1(\overline{X}^{an},j_*^m{\cal O}^*_{X^{an}})\not\simeq H^2_{DB}(X,\mathbb{Z}(1))$ in our case.\\

b) Let $X$ be the affine cubic surface in $\mathbb{C}^3$ defined by $z_1^3+z_2^3+z_3^3=-1$. Then we may take $\overline{X}$ as the corresponding projective cubic surface in $\mathbb{P}_3$ defined by $z_0^3+z_1^3+z_2^3+z_3^3=0$.\\
By [Ha] V Prop. 4.8, p. 401, we have $Pic\,\overline{X}\simeq \mathbb{Z}^7$. According to [Ha] II Prop. 6.5, p. 133, we have an exact sequence $\mathbb{Z}\to Pic\,\overline{X}\to Pic\,X\to 0$.\\ The hyperplane at infinity defines a Cartier divisor on $\overline{X}$ which is not principal, so $\mathbb{Z}\to Pic\,\overline{X}$ is not the zero map, hence injective. In fact, let $H\subset\mathbb{P}_3$ be the hyperplane at infinity. Then the corresponding first Chern class is $\neq 0$, and $H^2(\mathbb{P}_3^{an};\mathbb{Z})\to H^2(\overline{X}^{an};\mathbb{Z})$ is surjective, by the Lefschetz theorem on hyperplane sections.\\
Therefore we have a short exact sequence: $0\to \mathbb{Z}\to Pic\,\overline{X}\to Pic\,X\to 0$.\\
This implies that $rk\,Pic\,X=6$.\\
On the other hand, look at the arithmetic genus $p_a$ of $\overline{X}$: By [Ha] III Exc. 5.5, p. 231, we have $p_a=\dim\,H^2(\overline{X},{\cal O}_{\overline{X}})$, $H^1(\overline{X},{\cal O}_{\overline{X}})=0$, and by [Ha] I Exc. 7.2, p. 54, we have $p_a=\binom{2}{3}=0$.\\
Also, $X^{an}$ is the Milnor fibre of $z_1^3+z_2^3+z_3^3$, so $H^2(X^{an};\mathbb{Z})\simeq\mathbb{Z}^8$. \\
The exact sequence\\
$H^1(\overline{X},{\cal O}_{\overline{X}})\to H_{DB}^2(X,\mathbb{Z}(1))\to H^2(X^{an};\mathbb{Z})\to H^2(\overline{X},{\cal O}_{\overline{X}})$ \\
implies therefore that $H^2_{DB}(X,\mathbb{Z}(1))\simeq \mathbb{Z}^8\not\simeq Pic\,X$.\\
More precisely, $Pic\,X\simeq \mathbb{Z}^6$: We have an exact sequence\\
$0\to Pic\,X\to H^2_{DB}(X,\mathbb{Z}(1))\to \mathbb{Z}^2\to 0$.\\
This comes from the exact sequence of Proposition 2.6b): we have $H^3(\overline{X}^{an},X^{an};\mathbb{Z})\simeq H_1(D;\mathbb{Z})\simeq \mathbb{Z}^2$ 
because $D$ is an elliptic curve,
and the exponential sequence yields an exact sequence $H^2(\overline{X}^{an},{\cal O}_{\overline{X}^{an}})\to H^2(\overline{X}^{an},{\cal O}^*_{\overline{X}^{an}})\to H^3(\overline{X}^{an};\mathbb{Z})$.\\
Now $p_a=0$, see above, so $H^2(\overline{X}^{an},{\cal O}_{\overline{X}^{an}})=0$, and $H^3(\overline{X}^{an};\mathbb{Z})\simeq H_1(\overline{X}^{an};\mathbb{Z})=0$ because of Lefschetz theorem on hyperplane sections, hence 
$H^2(\overline{X}^{an},{\cal O}^*_{\overline{X}^{an}})=0$.\\
Moreover, $D$ is irreducible, hence $H^1(\overline{X}^{an},j_*^m{\cal O}^*_{X^{an}})\simeq H^2_{DB}(X,\mathbb{Z}(1))$; see Lemma 2.8.\\
In particular, $Pic\,X\not\simeq H^1(\overline{X}^{an},j_*^m{\cal O}^*_{X^{an}})$.\\

So both examples show that $Pic\,X$ and $H^2_{DB}(X,\mathbb{Z}(1))$ are different, in general. Example a) shows that $H^1(\overline{X}^{an},j_*^m{\cal O}^*_{X^{an}})$ and $H^2_{DB}(X,\mathbb{Z}(1))$ are in general different, too, whereas Example b) shows that the same holds for $Pic\,X$ and $H^1(\overline{X}^{an},j_*^m{\cal O}^*_{X^{an}})$.\\

c)  Let $k\ge 2$ and $X$ be a smooth hypersurface of degree $d\ge k+2$ in $\mathbb{P}_{k+1}$. Then, by [Ha] III Exc. 5.5, p. 231: $p_a=\dim\,H^k(X^{an},{\cal O}_{X^{an}})$, and by [Ha] I Exc. 7.2, p. 54: $p_a=\binom{d-1}{k+1}\neq 0$. On the other hand, in contrast to $H^k(X^{an},{\cal O}_{X^{an}})$, the group $H^k(X^{an};\mathbb{Z})$ is finitely generated. The exponential sequence implies that $H^k(X^{an},{\cal O}^*_{X^{an}})\neq 0$, whereas $H^k(X,{\cal O}^*_X)=0$ according to the Remark 2.2. So we have no GAGA principle for $H^k(X,{\cal O}^*_X)$, $k\ge 2$.\\
The situation is not better if we use \v{C}ech cohomology instead of ordinary (flabby) cohomology: Note that it is not clear whether these cohomology theories agree for ${\cal O}^*_X$, $X$ algebraic variety, because this sheaf is not coherent algebraic and $X$ is not paracompact. Anyhow, $\check{H}^2(X,{\cal O}^*_X)\to H^2(X,{\cal O}^*_X)$ is injective, see [Go], so for $k=2$ we obtain $\check{H}^k(X,{\cal O}^*_X)=0$, so there is no GAGA principle for $\check{H}^2(X,{\cal O}^*_X)$, too.\\

{\bf References:}\\

[B1] A. Beilinson: Higher regulators and values of $L$-functions [Original Russian]. J. Soviet Math. {\bf 30}, no. 2, 2036-2070 (1985).\\

[B2] A. Beilinson: Notes on absolute Hodge cohomology. In: Appl. of alg. $K$-theory to alg. geometry and number theory, Part I, pp. 35-68 (Boulder, Colo. 1983). Contemp. Math. {\bf 55} (1986).\\

[D] P. Deligne: Equations diff\'erentielles \`a points singuliers r\'eguliers. Lecture Notes in Math. {\bf 163}, Springer-Verlag: Berlin 1970.\\

[Do] A. Dold: Lectures on algebraic topology. Springer-Verlag: Berlin 1972.\\

[E] F. El Zein: Introduction \`a la th\'eorie de Hodge mixte. Hermann: Paris 1991.\\

[EV] H. Esnault, E. Viehweg: Deligne-Beilinson Cohomology. In: Beilinson's Conjectures on Special Values of $L$-Functions, ed. M. Rapoport, N. Schappacher, P. Schneider, pp. 43-91. Perspectives in Mathematics {\bf 4}. Academic Press: Boston, Mass. 1988.\\

[G] A. Grothendieck: On the de Rham cohomology of algebraic varieties. Publ. Math. IHES {\bf 29}, 95-103 (1966).\\

[Ga] P. Gajer: Geometry of Deligne cohomology. Invent. Math. {\bf 127}, 155-207 (1997).\\

[Go] R. Godement: Topologie alg\'ebrique et th\'eorie des faisceaux. Hermann: Paris 1958.\\

[GM] S. I. Gelfand, Yu. I. Manin: Methods of homological algebra. Springer-Verlag: Berlin 1996.\\

[Ha] R. Hartshorne: Algebraic Geometry. Springer-Verlag: New York 1977.\\

[HL1] H. A. Hamm, L\^e D. T.: On the Picard group for non-complete algebraic varieties. In: Singularit\'es Franco-Japonaises, ed. J.-P. Brasselet, T. Suwa, pp. 71-86. S\'em. et Congr\`es {\bf 10}. Soc. Math. France, 2005.\\

[HL2] H. A. Hamm, L\^e D. T.: Picard groups for line bundles with connection [in preparation]\\

[M] B. Malgrange: Regular connections after Deligne. In: A. Borel et al.: Algebraic $D$-modules, pp. 151-172. Academic Press: Boston 1987.\\

[N] M. Nagata: Embedding of an abstract variety into a complete variety. J. Math. Kyoto Univ. {\bf 2}, 1-10 (1962).\\

[S] M. Saito: Hodge conjecture and mixed motives. I. Complex geometry and Lie theory (Sundance, UT, 1989), 283–303, Proc. Sympos. Pure Math. {\bf53}, Amer. Math. Soc., Providence, RI, 1991.\\

[Sp] E. H. Spanier: Algebraic topology. McGraw-Hill: N.Y. 1966.\\

\footnotesize
Helmut A. Hamm\\
Mathematisches Institut der WWU\\
Einsteinstr. 62\\
48149 M\"unster\\
Germany\\
hamm@uni-muenster.de

\end{document}